\documentclass[10pt]{amsart}

\newtheorem{defin}{Definition}[subsection]
\usepackage[finalnew]{trackchanges}

\addeditor{Review1}
\addeditor{Review2}
\addeditor{Review3}
\addeditor{Review4}
\addeditor{rephrasing}

\usepackage[margin=1.2in]{geometry}
\usepackage{xcolor}

\theoremstyle{definition}
\newtheorem{thm}{Theorem}[subsection]
\theoremstyle{remark}
\usepackage{cite}
\numberwithin{equation}{section}


\title{A review on BGK models for gas mixtures of mono and polyatomic molecules}

\author{Marlies Pirner $^{1,\dagger,\ddagger}$
}

\address{%
$^{1}$ \quad Institute of Mathematics, University of Würzburg; marlies.pirner@mathematik.uni-wuerzburg.de
}

\begin{document}
\maketitle
\textbf{abstract: } We consider the socalled Bathnagar-Gross-Krook (BGK) model, an approximation of the Boltzmann equation, describing the time evolution of a single momoatomic rarefied gas and satisfying the same two main properties (conservation properties and entropy inequality). However, in practical applications one often has to deal with twoo additional physical issues. First,  a gas often does not consist of only one species, but it consists of mixture a mixture of different species. Second, the particles can store energy not only in translational degrees of freedom but also in internal degrees of freedom like rotations or vibrations (polyatomic molecules). Therefore, here, we will present recent BGK models for gas mixtures for mono- and polyatomic particles and the existing mathematical theory for these models.

\textbf{keywords:} multi-fluid mixture, kinetic model, BGK approximation, degrees of freedom in internal energy, existence of solutions, large-time behaviour


\section{Introduction}



\change[rephrasing]{ In this paper we shall concern ourselves with a kinetic description of gas mixtures. In the case of mono atomic molecules and two species this is traditionally done via the Boltzmann equation for the density distributions}{ In this paper we will concern ourselves with a kinetic description of gas mixtures. In the case of mono atomic molecules and two species this is usually done with the Boltzmann equation for the distribution functions} $f_1=f_1(x,v,t), f_2=f_2(x,v,t)$, see for example \cite{Cercignani, Cercignani_1975}. Here, $x \in \mathbb{R}^d$ and $v \in \mathbb{R}^d$ are the phase space variables, position and velocity of the particles, and $t \geq 0$ denotes the time. The Boltzmann equation for gas mixtures is of the form
 \begin{align*}
     \partial_t f_1 + v \cdot \nabla_x f_1 = Q_{11}(f_1,f_1)+ Q_{12}(f_1,f_2), \\
  \partial_t f_2 + v \cdot \nabla_x f_2 = Q_{22}(f_2,f_2)+ Q_{21}(f_2,f_1)
 \end{align*}
 where the collision operators $Q_{11}(f_1,f_1)$ and $Q_{22}(f_2,f_2)$ satisfy
\begin{align}
 \int Q_{kk}(f_k,f_k) \begin{pmatrix} 1 \\ m_k v \\ m_k |v|^2 \end{pmatrix} dv = 0 
\label{consone}
\end{align}
for $ k = 1,2$, and $Q_{12}(f_1,f_2)$ and $Q_{21}(f_2,f_1)$ satisfy
\begin{align}
\begin{split}
 \int Q_{12}(f_1,f_2) dv = 0, \quad \int Q_{21}(f_2,f_1) dv=0, \\
 \int ( m_1 v Q_{12}(f_1,f_2)+ m_2 v Q_{21}(f_2,f_1)) dv = 0, \\
  \int (m_1 |v|^2 Q_{12}(f_1,f_2)+m_2 |v|^2 Q_{21}(f_2,f_1) )dv = 0. 
\label{cons}
\end{split}
\end{align}
The properties \eqref{consone} ensure conservation of the number of particles, momentum and energy in interactions of one species with itself. The properties \eqref{cons} ensure conservation of the number of particles, total momentum and total energy in interactions of one species with the other species. For the proof see for example \cite{Sirovich}. 
In addition, the collision operator satisfy the inequalities \cite{Aoki}
\begin{align}
    \int Q_{kk}(f_k,f_k) \ln f_k dv \leq 0, \quad k=1,2 \label{h}
\\
    \int Q_{12}(f_1,f_2) \ln f_1 dv + \int Q_{21}(f_2,f_1) \ln f_2 dv \leq 0
\label{H}
\end{align}
In the first inequality \eqref{h} we have equality if and only if $f_k$ is equal to a Maxwell distribution given by 
\begin{align} 
M_k &= \frac{n_k}{\sqrt{2 \pi \frac{T_k}{m_k}}^d }  \exp \left({- \frac{|v-u_k|^2}{2 \frac{T_k}{m_k}}} \right), \quad k=1,2.
\label{M}
\end{align}
\add[Review1]{Note that in this paper we shall write $T_k$ instead of $k_B T_k$, where $k_B$ is Boltzmann's constant.} Here, for any $f_1, f_2: \Lambda \subset \mathbb{R}^d \times \mathbb{R}^d \times \mathbb{R}^+_0 \rightarrow \mathbb{R}$ with $(1+|v|^2)f_1,(1+|v|^2)f_2 \in L^1(\mathbb{R}^3), f_1,f_2 \geq 0$\textcolor{black}{,} we relate the distribution functions to  macroscopic quantities by mean-values of $f_k$,
\begin{align}
\int f_k(v) \begin{pmatrix}
1 \\ v  \\ m_k |v-u_k|^2 
\end{pmatrix} 
dv =: \begin{pmatrix}
n_k \\ n_k u_k \\ d n_k T_k
\end{pmatrix} , \quad k=1,2,
\label{moments}
\end{align} 
where $n_k$ is the number density, $u_k$ the mean velocity and $T_k$ the mean temperature of species \textcolor{black}{ $k$ ($k=1,2$).} \remove[Review1]{Note that in this paper we shall write $T_k$ instead of $k_B T_k$, where $k_B$ is Boltzmann's constant.}
Additionally to a Maxwell distribution in the second inequality \eqref{H} we have equality if and only if $u_1=u_2$ and $T_1=T_2$ \cite{Aoki}. In the following we refer to these inequalities including the characterization of equality as entropy inequalities of H-Theorem. 

 \add[Review3]{If we are close to equilibrium} \cite{PirnerD,Boscarino}, \change[rephrasing]{the complicated interaction terms of the Boltzmann equation can be simplified by a so called BGK approximation, consisting of a collision frequency $\nu_{ij} n_j$ multiplied by the deviation of the distributions from local Maxwell distribution.}{the complicated interaction terms of the Boltzmann equation can be simplified by a BGK approximation. This consists of a collision frequency $\nu_{ij} n_j$ multiplied by the deviation of the distribution functions from a local Maxwell distribution.}
 \begin{align} \begin{split} \label{BGK}
\partial_t f_1 + v \cdot \nabla_x  f_1   &= \nu_{11} n_1 (M_1-f_1) + \nu_{12} n_2 (M_{12}- f_1),
\\ 
\partial_t f_2 + v \cdot \nabla_x f_2 &=\nu_{22} n_2 (M_2 - f_2) + \nu_{21} n_1 (M_{21}- f_2), 
\end{split}
\end{align}
\add[Review3]{The collision frequencies per density $\nu_{kj}$ are assumed to be dependent only on $x$ and $t$ and not on the microscopic velocity $v$. For references taking into account also a dependency on the microscopic velocity $v$ see} \cite{Struchtrupp} for the one species case and \cite{Pirner_velocity} for the gas mixture case. Here, $M_k$ is given by \eqref{M} and $M_{12}, M_{21}$ are given by
\begin{align} 
\begin{split}
M_{12} &= \frac{n_{1}}{\sqrt{2 \pi \frac{T_{12}}{m_1}}^d }  \exp({- \frac{|v-u_{12}|^2}{2 \frac{T_{12}}{m_1}}}), 
\\
M_{21} &= \frac{n_{2}}{\sqrt{2 \pi \frac{T_{21}}{m_2}}^d }  \exp({- \frac{|v-u_{21}|^2}{2 \frac{T_{21}}{m_2}}}),
\end{split}
\label{BGKmix}
\end{align}
for suitable $ u_{12}, u_{21}, T_{12}, T_{21}$ such that the conservation properties \eqref{cons} are satisfied. In the literature, there is another type of approximation of the form 
\begin{align} \begin{split} \label{BGK2}
\partial_t f_k + v \cdot \nabla_x  f_k   &= \nu_k (M^{(k)} - f_k), \quad \nu_k=\sum_{j=1}^2 \nu_{kj} n_k 
\end{split}
\end{align}
with only one relaxation term taking into account both type of interactions, inter- and intra-species interactions, with one Maxwell distribution
\begin{align}
M^{(k)}(x,v,t) &= \frac{n_k}{\sqrt{2 \pi \frac{T^{(k)}}{m_k}}^d }  \exp({- \frac{|v-u^{(k)}|^2}{2 \frac{T^{(k)}}{m_k}}}),
\label{BGKmix2}
\end{align}
for suitable $n^{(k)}, u^{(k)}, T^{(k)}$ such that the conservation properties \eqref{cons} are satisfied.
\change[rephrasing]{BGK  models give rise to efficient numerical computations, which are asymptotic preserving, that is they remain efficient even approaching the hydrodynamic regime}{BGK  models give rise to efficient numerical computation, which are asymptotic preserving. This means that  they remain efficient even approaching the hydrodynamic regime} \cite{Puppo_2007, Jin_2010,Dimarco_2014, Bennoune_2008,  Bernard_2015, Crestetto_2012}. \add[Review1]{However, the BGK approximation has one drawback. This  is its incapability of reproducing the correct Boltzmann hydrodynamic regime in the asymptotic continuum limit. Therefore, a modified version called the ES-BGK approximation was suggested  by Holway for one species} \cite{Holway}. \add[Review1]{Then the H-Theorem of this model  was shown in} \cite{Perthame} \add[Review1]{and existence and uniqueness of  mild solutions in} \cite{Yun_mild}. \add[Review1]{Alternatively, the Shakov model} \cite{Shakov} \add[Review1]{and a BGK model with velocity dependent collision frequency} \cite{Struchtrupp}\add[Review1]{was suggested to achieve the correct Prandtl number. For the BGK model with velocity dependent collision frequency, it is shown that a power law for the collision frequency can also lead to the proper Prandtl number.  This BGK model with velocity dependent collision frequency should be also constructed in a way such that it satisfies the conservation properties. That this works one has to replace the Maxwell distribution by a different function, for details see} \cite{Struchtrupp}. \add[Review1]{For this model also an H-Theorem can be proven. The existence of these modified functions is proven in} \cite{Pirner_velocity}. \add[Review1]{However, since BGK models can be the basis to build extended models as ES-BGK models, Shakov models and BGK models with velocity dependent collision frequency, we will mainly review BGK models for gas mixtures in this paper}.

\change[rephrasing]{We are interested in extensions  of a BGK model to gas mixtures since in applications one often has to deal with mixtures instead of a single gas. Moreover, evolution of a polyatomic gas is very important in applications, for instance air consists of a gas mixture of polyatomic molecules. But most kinetic models modelling air deal with the case of a mono atomic  gas consisting of only one species. Therefore, in this paper, we give an overview over existing BGK models for gas mixtures in the mono- and polyatomic case and theoretical results for them.}{
In this paper, we are interested in extensions  of a BGK model to gas mixtures since in applications one often has to deal with mixtures instead of a single gas. Moreover, evolution of a polyatomic gas is very important in applications, for instance air consists of a gas mixture of polyatomic molecules. But most kinetic models modelling air deal with the case of a mono atomic  gas consisting of only one species.}


The outline of the paper is as follows: In section \ref{sec2} we will present typical ansatzes of modelling gas mixtures with the BGK model and a review on recent results concerning existence of solutions and large-time behaviour. In section \ref{sec3} we will give a summary over existing BGK models for gas mixtures of polyatomic molecules and a review on existing theoretical results.
 
\section{BGK models for gas mixtures}
\label{sec2}
In the following, we will present existing BGK models for gas mixtures. Then, we will give an overview over existing theoretical results (existence of solutions, large-time behaviour).
\subsection{Overview on existing BGK models for gas mixtures in the literature}
\change[rephrasing]{Here we shall focus on gas mixtures modelled via a BGK approach. In the literature one can find two types of models. Just like the Boltzmann equation for gas mixtures contains a sum of collision terms on the right-hand side, one type of BKG models also has a sum of  BGK-type interaction terms in the relaxation operator}{Here we will focus on gas mixtures modelled via a BGK approach. In the literature one can find two types of  BGK models for gas mixtures. Just like the Boltzmann equation for gas mixtures has a sum of collision terms on the right-hand side, one type of BKG models also contains a sum of  BGK-type relaxation terms on the right-hand side} \eqref{BGK}. \annote[Review2]{Examples}{References in ascending order} are the models of  Asinari \cite{asinari},   Cercignani \cite{Cercignani_1975}, Garzo, Santos, Brey \cite{Garzo1989}, Greene \cite{Greene}, Gross and Krook \cite{gross_krook1956}, Hamel \cite{hamel1965},  Sofena \cite{Sofonea2001},  and recent models by Bobylev, Bisi, Groppi, Spiga, Potapenko \cite{Bobylev}; Haack, Hauck, Murillo \cite{haack} and by Klingenberg, Pirner, Puppo \cite{Pirner}. The other type of models contains only one collision term on the right-hand side \eqref{BGK2}. Examples for this are Andries, Aoki and Perthame \cite{AndriesAokiPerthame2002}, and the models in \cite{Brull_2012, Groppi}.
\add[Review2]{A comparison of these models concerning their hydrodynamic limit can be found in} \cite{Boscarino}
\add[Review1]{There are also many results concerning the hydrodynamic limit via the Chapman Enskog expansion, see for example} \cite{AndriesAokiPerthame2002,Bernard_2015,Brull_2012,Crestetto_2012,haack} \add[Review1]{ and extensions to ES-BGK models, Shakov models and BGK models with velocity dependent collision frequency} \cite{Groppi,Brull,Blaga2,Pirner_velocity}.
\subsubsection{BGK models for gas mixtures with one collision term}
BGK models for gas mixtures with one interaction term \cite{AndriesAokiPerthame2002,Brull_2012,Groppi} on the right-hand side have the form \eqref{BGK2}.
Now, the interspecies velocities $u^{(k)}$ and temperatures $T^{(k)}$ in \eqref{BGKmix2} have to be determined such that the conservation of total momentum and total energy
\begin{align*}
\sum_{k}^2 \int \nu_k  (M^{(k)} - f_k) \begin{pmatrix} v \\ |v|^2 \end{pmatrix} dv = 0
\end{align*}
as in \eqref{cons} is satisfied. 
This gives $d+1$ constraints for the $2 (d+1)$ quantities $u^{(k)}$, $T^{(k)}$, $k=1,2$. Therefore there is additional freedom to choose these quantities. Examples in the literature are given in \cite{AndriesAokiPerthame2002,Brull_2012,Groppi}. In the first case \cite{AndriesAokiPerthame2002} the quantities $u^{(k)}$ and $T^{(k)}$ are chosen such that the exchange terms of momentum and energy 
\begin{align*}
\int \nu_k  (M^{(k)} - f_k) \begin{pmatrix} v \\ |v|^2 \end{pmatrix} dv 
\end{align*}
coincide with the exchange terms of momentum and energy of the Boltzmann equation for Maxwell molecules. For the details, see \cite{AndriesAokiPerthame2002}. This leads to the choice
\begin{align}
\begin{split}
\label{choice12}
u^{(k)}&= u_k + \sum_{j=1}^2 2 \frac{\chi_{kj}}{\nu_k} \frac{m_j}{m_k+m_j} n_j (u_j-u_k)
\\
T^{(k)} &= T_k - \frac{m_k}{d} |u^{(k)}- u_k|^2 \\ &+ \sum_{j=1}^2 2 \frac{m_k m_j}{m_k+m_j} \frac{\chi_{kj}}{\nu_{kj}} n_j \frac{2}{m_k+m_j} (T_j-T_k + \frac{m_j}{d} |u_j-u_k|^2)
\end{split}
\end{align}\note[Review1]{Typo N2 changed into sum until 2}
where $\chi_{12},$ $\chi_{21}$ are parameters which are related to the differential cross section. For the detailed expressions see \cite{AndriesAokiPerthame2002}.

The model also satisfies the conservation properties \eqref{cons} and the H-theorem \eqref{H}
\change[rephrasing]{ with equality if and only if the distribution functions are Maxwell distributions with equal mean velocity and temperature. }
{. In the H-theorem one has equality if and only if the distribution functions are Maxwell distributions with the same mean velocity and temperature. }
The model of Andries, Aoki and Perthame  has another property (see proposition 3.2 in \cite{AndriesAokiPerthame2002}). 
It is called the indifferentiability principle. \change[rephrasing]{It denotes the following property: When the masses $m_k, ~ k=1,2$ and the collision frequencies $\nu_{kj}, ~ k,j=1,2$ are  identical, the total distribution function $f=f_1+f_2$ obeys a single species BGK equation.}{This means the following. When the masses $m_k, ~ k=1,2$ and the collision frequencies $\nu_{kj}, ~ k,j=1,2$ are  the same for each species, the total distribution function $f=f_1+f_2$ satisfies a single species BGK equation.
}

\add[Review2]{A derivation of the Navier-Stokes system in the compressible regime and the corresponding transport coefficients can be found in section 4 of} \cite{AndriesAokiPerthame2002}.




Another model in the literature with shape \eqref{BGK2} is the model \cite{Groppi}. Here $u^{(k)}, T^{(k)}$ are chosen such that all interspecies velocities $u^{(k)}$ and temperatures $T^{(k)}$ are equal
\begin{align*}
u^{(k)} = \bar{u}, \quad T^{(k)} = \bar{T}
\end{align*}
for all $k=1,2$, where $\bar{u}$ and $\bar{T}$ are determined such that the conservation properties as in \eqref{cons} are satisfied. This leads to the choice
\begin{align*}
\bar{u}= \frac{\sum_{s=1}^2 \nu_s m_s n_s u_s}{\sum_{s=1}^2 \nu_s m_s n_s}, \quad \bar{T}= \frac{ \sum_{s=1}^2 \nu_s n_s (m_s (|u_s|^2 - |\bar{u}|^2) + d T_s)}{d \sum_{s=1}^2 \nu_s n_s}
\end{align*}
In \cite{Groppi2}, it is proved that the positivity of the temperature $\bar{T}$ is guaranteed and the H-Theorem holds for the space-homogeneous case. 
\add[Review2]{The hydrodynamic limit and corresponding transport coefficients of this models can be found in section 5 in} \cite{Boscarino}.

Another model with shape \eqref{BGK2} is the model \cite{Brull_2012}. Here the aim was to derive the BGK model for gas mixtures from an entropy minimization principle ensuring that the model satisfies the exact Fick and Newtons laws in the hydrodynamical limit to the Navier-Stokes equations. This leads to a choice of different values for $u^{(k)}$ and equal values $T^{(k)} = T^*$ for all $k=1,2$ for the temperatures. For the detailed expressions, see \cite{Brull_2012}. \add[Review2]{The transport equations of the hydrodynamic regime for this model can be found in section 5 of} \cite{Brull_2012}.
\subsubsection{BGK models for gas mixtures with two collision terms}
Now, we review BGK models for gas mixtures with two collision terms.
%
%
%
In case of a gas mixture, \change[rephrasing]{the particles of one species can interact with either themselves or with  particles of the other species. One takes this into account  by writing two interaction terms in the equations}{if we assume that we only have binary interactions, there are two possibilities. The particles of one species can interact with themselves or with  particles of the other species. One can take this into account  by writing two interaction terms in the equations} \eqref{BGK}. \change[rephrasing]{This means that the right-hand side of the equations  consists of a sum of two collision operators. This structure is also described in}{This means that the right-hand side of the equations now consists of a sum of two relaxation operators. This structure is also described in} \cite{Cercignani, Cercignani_1975}.
\change[rephrasing]{ This leads  to  two types of equilibrium distributions. Due to the interaction of a species $k$ with itself, we expect a relaxation towards an equilibrium distribution $M_k$. And due to the interaction of a species with another one, we expect a relaxation towards a different equilibrium distribution $M_{kj}$.}{ This leads  to  two different types of equilibrium distributions. Due to an interaction of a species $k$ with itself, we expect a relaxation to an equilibrium distribution $M_k$. And due to the interaction of a species with the other species, we expect a relaxation towards a different mixture equilibrium distribution $M_{kj}$. }

%
\change[rephrasing]{The quantities $\nu_{kk} n_k$  are the collision frequencies of the particles of each species with itself, while $\nu_{kj} n_j$ are related to interspecies collisions.
To be flexible in choosing the relationship between the collision frequencies, we now assume the relationship}{The quantities $\nu_{kk} n_k$  are the one species collision frequencies, while the collision frequencies $\nu_{kj} n_j$ are related to interspecies interactions.
To be flexible in choosing the relationship between the interspecies collision frequencies, we assume the following relationship} 
\begin{equation} 
\nu_{12}=\varepsilon \nu_{21}, \quad 0 < \varepsilon \leq 1.
\label{coll}
\end{equation}
\change[rephrasing]{The restriction on $\varepsilon$ is without loss of generality. If $\varepsilon >1$, exchange the notation $1$ and $2$ and choose $\frac{1}{\varepsilon}.$} {The restriction on $\varepsilon$ is without loss of generality. If $\varepsilon >1$, we can exchange the notation $1$ and $2$ and choose $\frac{1}{\varepsilon}$ instead.} 

\add[Review2]{Let us provide an example. We consider a plasma with electrons and ions. Let us first denote the electrons with the index $e$ and ions with the index $i$. Then a common relationship found in the literature} \cite{bellan2006} \add[Review2]{is $\nu_{ie}= \frac{m_e}{m_i} \nu_{ei}$ or equivalent $\nu_{ei}=\frac{m_i}{m_e} \nu_{ie}$. Now, if we want to use the notation $1$ and $2$, we have two possibilities. The first one is to choose the notation 1 for electrons and 2 for ions. In this case the mass ratio of the two particles is $\frac{m_2}{m_1} >>1$, and we get $\nu_{12}= \frac{m_2}{m_1} \nu_{21}$. So we have $\varepsilon = \frac{m_2}{m_1}>1$. The other possibility is to choose the notation 1 for ions and 2 for electrons. In this case the mass ratio of the two kinds of particles is $\frac{m_2}{m_1} <<1,$ and we have $\nu_{12}= \frac{m_2}{m_1} \nu_{21}$, therefore $\varepsilon=\frac{m_2}{m_1}<1$. So in this case, we would use the second choice for the notation. The condition} \eqref{coll} \add[Review2]{ will enter in the proof of the H-Theorem. }

 \change[rephrasing]{In addition, we assume that all collision frequencies are positive.
The Maxwell distributions $M_k$ in}{ Additionally, we assume that all collision frequencies are strictly positive. 
The Maxwell distribution $M_k$ in} \eqref{BGKmix}
 \change[rephrasing]{ have the same density, mean velocity and temperature as $f_k$, respectively. With this choice, we guarantee the conservation of the number of particles, momentum and energy in interactions of one species with itself (see section 2.2 in}{has the same density, mean velocity and temperature as $f_k$. With this choice, it can be guaranteed that we have conservation of the number of particles, momentum and energy in interactions of a species with itself (see section 2.2 in} \cite{Pirner}).
\change[rephrasing]{The remaining parameters $ u_{12}, u_{21}, T_{12}, T_{21}$  will be determined using conservation of total momentum and energy. 
Due to the choice of the densities, we have conservation of the number of particles, see Theorem 2.1 in}{The remaining parameters $ u_{12}, u_{21}, T_{12}, T_{21}$  will now be determined using conservation of total momentum and total energy. 
Due to the choice of the densities, one can prove conservation of the number of particles, see Theorem 2.1 in} \cite{Pirner}.
\change[rephrasing]{If we further assume that $u_{12}$ is a linear combination of $u_1$ and $u_2$}{We further assume that $u_{12}$ is a linear combination of $u_1$ and $u_2$}
 \begin{align}
u_{12}= \delta u_1 + (1- \delta) u_2, \quad \delta \in \mathbb{R},
\label{convexvel}
\end{align} then we have conservation of total momentum 
provided that
\begin{align}
u_{21}=u_2 - \frac{m_1}{m_2} \varepsilon (1- \delta ) (u_2 - u_1),
\label{veloc}
\end{align}
see Theorem 2.2 in \cite{Pirner}.
If we further assume that $T_{12}$ is of the following form
\begin{align}
\begin{split}
T_{12} &=  \alpha T_1 + ( 1 - \alpha) T_2 + \gamma |u_1 - u_2 | ^2,  \quad 0 \leq \alpha \leq 1, \gamma \geq 0 ,
\label{contemp}
\end{split}
\end{align}
then we have conservation of total energy 
provided that
\begin{align}
\begin{split}
T_{21} =\left[ \frac{1}{d} \varepsilon m_1 (1- \delta) \left( \frac{m_1}{m_2} \varepsilon ( \delta - 1) + \delta +1 \right) - \varepsilon \gamma \right] |u_1 - u_2|^2 \\+ \varepsilon ( 1 - \alpha ) T_1 + ( 1- \varepsilon ( 1 - \alpha)) T_2,
\label{temp}
\end{split}
\end{align}
see Theorem 2.3 in \cite{Pirner}.
In order to ensure the positivity of all temperatures, we need to restrict $\delta$ and $\gamma$ to 
 \begin{align}
0 \leq \gamma  \leq \frac{m_1}{d} (1-\delta) \left[(1 + \frac{m_1}{m_2} \varepsilon ) \delta + 1 - \frac{m_1}{m_2} \varepsilon \right],
 \label{gamma}
 \end{align}
and
\begin{align}
 \frac{ \frac{m_1}{m_2}\varepsilon - 1}{1+\frac{m_1}{m_2}\varepsilon} \leq  \delta \leq 1,
\label{gammapos}
\end{align}
see Theorem 2.5 in \cite{Pirner}. For this model, one can prove an H-theorem as in \eqref{H} with equality if and only if $f_k, ~k=1,2$ are Maxwell distributions with equal mean velocity and temperature, see \cite{Pirner}.

This model contains a lot of proposed models in the literature as special cases.
\annote[Review2]{Examples}{References in ascending order} are the models of  Asinari \cite{asinari},   Cercignani \cite{Cercignani_1975}, Garzo, Santos, Brey \cite{Garzo1989}, Greene \cite{Greene}, Gross and Krook \cite{gross_krook1956}, Hamel \cite{hamel1965},  Sofena \cite{Sofonea2001},  and recent models by Bobylev, Bisi, Groppi, Spiga, Potapenko \cite{Bobylev}; Haack, Hauck, Murillo \cite{haack}.

The second last \cite{Bobylev} presents an additional motivation how it can be derived formally from the Boltzmann equation. The last one \cite{haack} presents a Chapman-Enskog expansion \add[Review2]{ with transport coefficients in section 5, a comparison with other BGK models for gas mixtures in section 6  and a numerical implementation in section 7.} 
\subsection{Theoretical results of BGK models for gas mixtures}
In this section, we present theoretical results for the models presented in section 2.1. \add[Review3]{We start with reviewing some existing theoretical results for the one-species BGK model.} \add[Review4]{Concerning the existence of solutions the first result was proven by Perthame  in} \cite{Perth}. \add[Review4]{It is a result on global weak solutions for general initial data. This result was inspired by Diperna and Lion from a result on the Boltzmann equation} \cite{DL}. 
In \cite{Perthame}, \add[Review4]{the authors consider mild solutions and also  obtain the uniqueness in the periodic bounded domain. There are also results of stationary solutions on a 1-dimensional finite interval with inflow boundary conditions in} \cite{Ukai}. \add[Review4]{In a regime near a global Maxwell distribution, the global existence in the whole space $\mathbb{R}^3$ was established in} \cite{Yun1}. 
\add[Review3]{Concerning convergence to equilibrium Desvillettes proved strong convergence to equilibrium considering the thermalizing effect of the wall for reverse and specular reflection boundary condition in a periodic box} \cite{Des}. In  \cite{Saint}, \add[Review4]{the fluid limit of the BGK model is considered.} 

In the following we will present theoretical results for BGK models for gas mixtures.
\subsubsection{Existence of solutions}
First, we will present an existence result of mild solutions under the following assumptions for both type of models.
\begin{enumerate}
\item We assume periodic boundary conditions in $x$. Equivalently\textcolor{black}{,} we can construct solutions satisfying 
{\small
$$f_k(t,x_1,..., x_d, v_1,...,v_d)= f_k(t,x_1,...,x_{i-1},x_i + a_i,x_{i+1},...x_d,v_1,...v_d)$$}
for all $i=1,...,d$ and a suitable $\lbrace a_i\rbrace \in \mathbb{R}^d$ with positive components, for $k=1,2$.
\item We require that the initial values $f_k^0, i=1,2$ satisfy assumption $1$.
\item We are on the bounded domain in space $\Lambda=\lbrace x \in \mathbb{R}^N | x_i \in (0,a_i)\rbrace$.
\item Suppose that $f_k^0$ satisfies $f_k^0 \geq 0$, $(1+|v|^2) f_k^0 \in L^1(\Lambda \times \mathbb{R}^d)$ with \\$\int f_k^0 dx dv =1, k=1,2$.
\item Suppose $N_q(f_k^0):= \sup f_k^0(x,v)(1+|v|^q) = \frac{1}{2} A_0 < \infty$ for some $q>d+2$.
\item Suppose $\gamma_k(x,t):= \int f_k^0(x-vt,v) dv \geq C_0 >0$ for all $t\in\mathbb{R}.$
\item Assume that the collision frequencies are written as 
\begin{align}
\begin{split}
\nu_{jk}(x,t) n_k(x,t) = \tilde{\nu}_{jk}\frac{n_k(x,t)}{n_j(x,t)+n_k(x,t)}, \quad j,k=1,2,
\end{split}
\label{cond_coll}
\end{align}
with constants $\tilde{\nu}_{jk}>0$. 
\end{enumerate}
\label{ass}
With these assumptions\textcolor{black}{,} we can show the following Theorem, existence of mild solutions in the following sense.
\begin{defin}\label{milddef}\note[Review4]{I added this definition}
We call $(f_1, f_2)$ with $(1+|v|^2)f_k \in L^1(\mathbb{R}^N), f_1,f_2 \geq 0$ a mild solution to \eqref{BGK} under the conditions of the collision frequencies \eqref{cond_coll} iff $f_1,f_2$ satisfy
{\footnotesize
\begin{align}
\begin{split}
&f_k(x,v,t)= e^{-\alpha_k(x,v,t)} f_k^0(x-tv,v) \\ &+ e^{-\alpha_k(x,v,t)} \int_0^t [ \tilde{\nu}_{kk} \frac{n_k(x+(s-t)v,s)}{n_k(x+(s-t)v,s)+ n_j(x+(s-t)v,s)} M_k(x+(s-t)v,v,s) \\ &+ \tilde{\nu}_{kj} \frac{n_j(x+(s-t)v,s)}{n_k(x+(s-t)v,s)+ n_j(x+(s-t)v,s)} M_{kj}(x+(s-t)v,v,s)]] e^{\alpha_k(x+(s-t)v,v,s)} ds,
\end{split}
\end{align}}
where $\alpha_k$ is given by
{\small
\begin{align}
\begin{split}
\alpha_k(x,v,t) = \int_0^t [\tilde{\nu}_{kk} \frac{n_k(x+(s-t)v,s)}{n_k(x+(s-t)v,s) + n_j(x+(s-t)v,s) }\\ +\tilde{\nu}_{kj} \frac{n_j(x+(s-t)v,s)}{n_k(x+(s-t)v,s) + n_j(x+(s-t)v,s) } ] ds,
\end{split}
\end{align}}
for $k,j =1,2, ~k\neq j$.
\end{defin}

 The proof can be found in \cite{Pirner2}.The main idea consists in proving Lipschitz continuity of the Maxwell distribution $M_{kj}$ and bounds on the macroscopic quantities needed for this. 
\begin{thm}
Under the assumptions 1.-7., there exists a unique non-negative mild solution $(f_1,f_2)\in C(\mathbb{R}^+ ; L^1((1+ |v|^2) dv dx)$ of the initial value problem \eqref{BGK} with \eqref{moments}, \eqref{convexvel}, \eqref{veloc}, \eqref{contemp} and \eqref{temp}, and to the initial value problem to \eqref{BGK2} with \eqref{choice12} . Moreover, for all $t>0$ the following bounds hold:
\begin{align*}
|u_k(t)|, |u_{kj}(t)| \leq A(t) < \infty, \quad
n_k(t) \geq C_0 e^{-t} >0, \quad
T_k(t), T_{kj}(t) \geq B(t)>0,
\end{align*} 
for $k,j=1,2, k \neq j$ and some constants $A(t),B(t)$.
\end{thm}
\subsubsection{Large-time behaviour}
In this section, we will give an overview over existing results on the large- time behaviour for BGK models for gas mixtures. 
 We denote by $H(f)= \int f \ln f dv$ the entropy of a function $f$ and by $H(f|g)= \int f \ln \frac{f}{g} dv$ the relative entropy of $f$ and $g$. Then, one can prove the following theorems. The proofs are given in \cite{Crestetto_2012}. 
 
\begin{thm}
In the space homogeneous case for the model \eqref{BGK} with \eqref{moments}, \eqref{convexvel}, \eqref{veloc}, \eqref{contemp} and \eqref{temp} 
we have the following decay rate of the distribution functions $f_1$ and $f_2$
$$ || f_k - M_k ||_{L^1(dv)} \leq 4 e^{- \frac{1}{2} C t } [ H(f_1^0|M_1^0) + H(f_2^0 | M_2^0)]^{\frac{1}{2}},  \quad k=1,2$$
where $C$ is a constant given by 
$$ C= \min \lbrace \nu_{11} n_1 + \nu_{12} n_2,..., \nu_{21} n_1 + \nu_{22} n_2 \rbrace,$$
and the index $0$ denotes the value at time $t=0$.
\end{thm}
The main ingredient is to prove the inequality 
$$ \nu_{12} n_2 H(M_{12}) + \nu_{21} n_2 H(M_{21}) \leq \nu_{12} n_2 H(M_1) + \nu_{21} n_1 H(M_2)$$
Therefore, this theorem can also be proven in a similar way for the model \eqref{BGK2} with \eqref{choice12}, since  a corresponding inequality for the model \eqref{BGK2} with \eqref{choice12} of the form
\begin{align*}
\nu_1 H(M^{(1)}) + \nu_2 H(M^{(2)}) \leq \nu_1 H(M_1) + \nu_2 H(M_2)
\end{align*}
is proven in \cite{AndriesAokiPerthame2002}.
The next two theorems can also be easily extended to the model \eqref{BGK2} with \eqref{choice12} because it satisfies the same macroscopic behaviour as the model \eqref{BGK} with the choice 
\begin{align*}
\delta &= -2 \frac{m_2}{m_1+m_2} \frac{\chi_{12}}{\nu_{12}}+1, \\
\alpha &= - 4 \frac{m_1 m_2}{(m_1 + m_2)^2} \frac{\chi_{12}}{\nu_{12}} +1, \\
\gamma &= \frac{4}{3} \frac{m_1 m_2^2}{(m_1+m_2)^2} \frac{\chi_{12}}{\nu_{12}} n_1 n_2(1- \frac{\chi_{12}}{\nu_{12}}).
\end{align*}
  \begin{thm}\label{th:estimate_vel}
Suppose that $\nu_{12}$ is constant in time. In the space-homogeneous case of the model \eqref{BGK} with \eqref{moments}, \eqref{convexvel}, \eqref{veloc}, \eqref{contemp} and \eqref{temp}, we have the following decay rate of the velocities
\begin{equation*}
|u_1(t) - u_2(t)|^2 = e^{- 2 \nu_{12} (1- \delta)\left(n_2+\frac{m_1}{m_2} n_1\right) t} |u_1(0) - u_2(0)|^2. 
\end{equation*}
\end{thm}
  \begin{thm}\label{th:estimate_temp} 
Suppose $\nu_{12}$ is constant in time. In the space-homogeneous case of the model \eqref{BGK} with \eqref{moments}, \eqref{convexvel}, \eqref{veloc}, \eqref{contemp} and \eqref{temp}, we have the following decay rate of the temperatures
\begin{equation*}
\begin{split}
\textcolor{black}{T_1(t) - T_2(t) = 
 e^{- C_1t} \left[T_1(0) - T_2(0)+\frac{C_2}{C_1-C_3 } ( e^{(C_1-C_3) t} - 1) |u_1(0) - u_2(0)|^2 \right],} 
 \end{split}
\end{equation*}
where the constants are defined by
\begin{align*}
C_1&=(1- \alpha)\nu_{12}\left(n_2+n_1\right),\\
C_2&=\nu_{12}\left(n_2\left((1-\delta)^2+\frac{\gamma}{m_1}\right)-n_1\left(1-\delta^2-\frac{\gamma}{m_1}\right)\right),\\
C_3&=2 \nu_{12} (1- \delta)\left(n_2+\frac{m_1}{m_2} n_1\right).
\end{align*}
\end{thm}
The proofs can be found in \cite{Crestetto_2012}.
There are also results in the space-inhomogeneous case for linearized collision operator of the model \eqref{BGK} with \eqref{moments}, \eqref{convexvel}, \eqref{veloc}, \eqref{contemp} and \eqref{temp} for two species, see \cite{Pirner_Liu}.
For this, we consider a solution $(f_1,f_2)$ to \eqref{BGK} which is close to the equilibrium $(f_1^{\infty}, f_2^{\infty})$ with 
\begin{align}
f_k(x,v,t) = f_k^{\infty}(v) + h_k(x,v,t), \quad f_k^{\infty} (v) = \frac{n_{\infty,k}}{(2 \pi /m_k)^{d/2}} \exp \left(- \frac{|v|^2}{2 /m_k}\right)
\label{ansatz}
 \end{align}
Then, we have
\begin{align}
\begin{split}
 &n_k(x,t)= n_{\infty,k} + \sigma_k(x,t) \quad \text{with} \quad \sigma_k(x,t) = \int h_k(x,v,t) dv \\
&(n_k u_k)(x,t) = \int v f_k(x,v,t) dv = \mu_k(x,t) \quad \text{with} \quad \mu_k(x,t) = \int v h_k(x,v,t) dv \\ &P_k(x,t) = \frac{m_k}{d}  \int |v-u_k|^2 f_k(x,v,t) dv = n_{\infty,k} + \frac{1}{d}\left[\tau_k(x,t) - \frac{m_k |\mu_k(x,t)|^2}{n_{\infty,k}+ \sigma_k(x,t)}\right] \\ &\hspace{5cm}\quad \text{with} \quad \tau_k(x,t) = m_k \int |v|^2 h_k(x,v,t) dv.
\end{split}
\label{momh}
\end{align} 
Now, we do Taylor expansion of the terms $M_1, M_2, M_{12}, M_{21}$ with respect to $\sigma_1, \sigma_2, \mu_1, \mu_2, \tau_1$ and $\tau_2$ around 0 and only take first order terms. Moreover, one neglects quadratic terms of the form $\sigma_k \sigma_l, \sigma_k \mu_l$ and $\sigma_k \tau_l$. Then one obtains the linearized system
{\footnotesize
\begin{align} \begin{split} \label{BGKhsmall}
&\quad\partial_t h_1 + v \cdot \nabla_x  h_1 \\ &= \nu_{11} n_{\infty,1}\left(f_1^{\infty} (v) [(\frac{1+d/2}{n_{\infty,1}} - \frac{m_1 |v|^2}{2 n_{\infty,1}}) \sigma_1(x,t) + \frac{m_1}{n_{\infty,1}} v \cdot \mu_1(x,t) + \frac{1}{n_{\infty,1}} (- \frac{1}{2} + \frac{m_1|v|^2}{2d}) \tau_1(x,t)]- h_1\right) \\&\quad + \nu_{12} n_{\infty,2}  (f_1^{\infty} [ \frac{1}{n_{\infty,1}}(1+ \frac{\alpha}{2} (d- m_1 |v|^2)) \sigma_1 + \frac{1}{2} \frac{1}{n_{\infty,2}}(1-\alpha) (d- m_1 |v|^2) \sigma_2 \\&\quad + \frac{1}{n_{\infty,1}}\delta m_1 v \cdot \mu_1 + \frac{1}{n_{\infty,2}}(1-\delta) m_1 v \cdot \mu_2 + \frac{1}{2} \frac{1}{n_{\infty,1}}\alpha (\frac{1}{d} m_1 |v|^2 -1) \tau_1 +  \frac{1}{2} \frac{1}{n_{\infty,2}}(1-\alpha) (\frac{1}{d} m_1 |v|^2 -1) \tau_2] - h_1),
\\[10pt]
&\quad\partial_t h_2 + v \cdot \nabla_x h_2 \\&  =\nu_{22} n_{\infty,2}\left(f_2^{\infty} (v) [(\frac{1+d/2}{n_{\infty,2}} - \frac{m_2 |v|^2}{2 n_{\infty,2}}) \sigma_2(x,t) + \frac{m_2}{n_{\infty,2}} v \cdot \mu_2(x,t) + \frac{1}{n_{\infty,2}} (- \frac{1}{2} + \frac{m_2|v|^2}{2d}) \tau_2(x,t)] -h_2\right)\\&\quad  + \nu_{21} n_{\infty,1} (f_2^{\infty} [\frac{1}{2} \frac{1}{n_{\infty,1}}\varepsilon(1-\alpha) (d- m_2 |v|^2) \sigma_1+ \frac{1}{n_{\infty,2}}(1+ \frac{1-\varepsilon(1-\alpha)}{2} (d- m_2 |v|^2)) \sigma_2  
 \\&\quad  + \frac{1}{n_{\infty,1}}\varepsilon (1- \delta) m_1 v \cdot \mu_1 + \frac{1}{n_{\infty,2}}(1-\frac{m_1}{m_2} \varepsilon (1- \delta)) m_2 v \cdot \mu_2 + \frac{1}{2} \frac{1}{n_{\infty,1}}\varepsilon (1- \alpha) (\frac{1}{d} m_2 |v|^2 -1) \tau_1 \\ &\quad +  \frac{1}{2} \frac{1}{n_{\infty,2}}(1-\varepsilon (1- \alpha)) (\frac{1}{d} m_2 |v|^2 -1) \tau_2] -h_2). 
\end{split}
\end{align}}

and the following hypocoercivity result
\begin{thm}
Let $x$ in the $d$-dimensional torus of side length $L$. For each side length $L>0$ and dimension $d=1$, there exists an entropy functional $e(f_1, f_2)$  satisfying 
$$ c_d(L)\, e(f_1, f_2) \leq ||f_1 - f_1^{\infty}||^2_{L^2\left((\frac{f_1^{\infty}(v)}{n_{\infty,1}})^{-1}dv dx\right)}+ ||f_2 - f_2^{\infty}||^2_{L^2\left((\frac{f_2^{\infty}(v)}{n_{\infty,2}})^{-1}dv dx\right)} \leq C_d(L)\, e(f_1, f_2)$$
with some positive constants $c_d$, $C_d$ \textcolor{black}{that depend on $L$.}

Moreover, assume that $$\nu_{11} n_{\infty,1} + \nu_{12} n_{\infty,2}= 1 \qquad\text{and }\qquad \nu_{22} n_{\infty,2} + \nu_{21} n_{\infty,1}= 1, $$ then any solution $(h_1(t), h_2(t))$ to \eqref{BGKhsmall} in dimension $d=1$ with $e(h_1(0), h_2(0))< \infty,$ normalized according to
\begin{align}
\begin{split}
&\int \sigma_1(x,0) dx = \int \sigma_2(x,0) dx = 0, 
\\& \int (m_1 \mu_1(x,0) + m_2 \mu_2(x,0)) dx = 0, \,  \int (\tau_1(x,0) + \tau_2(x,0) ) dx = 0,
\end{split}
\label{cons2}
\end{align}
 then satisfies
{\small $$e(h_1(t), h_2(t)) \leq e^{-\tilde{C} t}
e(h_1(0), h_2(0)), $$}
where $\tilde{C}$ is given by {\scriptsize $$\tilde{C}=2 \min \lbrace\nu_{12} n_{\infty,2} (1- \delta), \nu_{12} n_{\infty,2} (1- \alpha), \nu_{11} n_{\infty,1}+ \nu_{12} n_{\infty,2}, \nu_{12} n_{\infty,1} \frac{m_1}{m_2} (1- \delta),  \nu_{12} n_{\infty,1} (1- \alpha), \nu_{22} n_{\infty,2} + \nu_{12} n_{\infty,1} , 2 \mu\rbrace.$$ }
Here, \textcolor{black}{$\mu^d(L)$} is a one species decay rate developed in \textcolor{black}{theorem 1.1 in } \cite{achleitner2017multi}. 
\label{theo}
\end{thm}
\section{BGK models for gas mixtures of polyatomic molecules}
\label{sec3}
In this section, we review recent models for gas mixtures of polyatomic molecules. This means that we take into account that the particles can also store energy in degrees of freedom in internal energy like rotations and vibrations. We start with an overview on the one species models to introduce typical modeling aspects in this context. 

First, we introduce a dependency on the degrees of freedom in internal energy. In the literature, this is done in different ways:
\begin{itemize}
\item Discrete dependency on the degrees of freedom in internal energy
\end{itemize}
Let us consider a system with $l>0$ internal energy states $E_l$. Then, we introduce $l$ distribution functions $f^{(1)}_l=f^{(1)}_l(x,v,t,E_l)$, one distribution function for each internal energy state $E_l$. This is for example considered in \cite{Morse}.
\begin{itemize}
\item A continuous scalar dependency on the degrees of freedom in internal energy
\end{itemize}
In this description we take into account degrees of freedom in internal energy by introducing an internal energy parameter $I$ which takes into account all degrees of freedom with a continuous scalar variable. Then, we introduce a distribution function $f^{(2)}=f^{(2)}(x,v,t, I)$. This is done for example in \cite{AndriesPerthame2001}.
\begin{itemize}
\item A discrete and continuous dependency of the degrees of freedom in internal energy
\end{itemize}
In this description, the rotations and vibrations of a diatomic gas are described in different ways. Here we define the distribution function $f^{(3)}=f^{(3)}(x,v,t, I,i)$ where $I \in \mathbb{R}^+_0$ describes the internal energy in a continuous way, whereas $i$ denotes the $i$th vibrational energy level of the corresponding vibrational energy  $\frac{ i R h \nu}{k_B}$ ($h$ is the Planck constant, $R$ the fundamental gas constant, while $\nu$ is the fundamental vibrational frequency of the molecule). This is done for example  in \cite{Mathiaud}.
\begin{itemize}
\item A vector-valued continuous dependency of the degrees of freedom in internal energy
\end{itemize}
Last, we consider a distribution function $f^{(4)}=f^{(4)}(x,v,t, \eta)$ where $\eta \in \mathbb{R}^l$ is a vector-valued continuous variable for the internal degrees of freedom, one component for each degree of freedom. Then $|\eta|^2$ has the meaning of microscopic energy in the internal degrees of freedom. Especially, one can also describe rotations and vibrations in a separate way. This is introduced in \cite{Bernard}. \\

In this article, we focus on the treatment of the additional continuous variable. For this we consider a distribution function $f(x,v,t, \mathcal{E})$ where $\mathcal{E}$ represents the dependency on the internal degrees of freedom and can be either a scalar $I$ or a vector $\eta$. So $f$ can be either $f^{(2)}$, $f^{(3)}$ for a fixed $i$ or $f^{(4)}$. Then we define the macroscopic quantities as

\begin{align}
\int f(v, \mathcal{E}) \begin{pmatrix}
1 \\ v \\ \mathcal{E} \\ m |v-u|^2 \\ m e(\mathcal{E}) 
\end{pmatrix} 
dv d\mathcal{E}=: \begin{pmatrix}
n \\ n u \\ n \bar{\mathcal{E}}\\ d  T_{tr} \\ l n T_{int} 
\end{pmatrix} , 
\label{momentsone}
\end{align} 
where we have $e(I)=I^{2/l}$ in the scalar case and $e(\eta)=|\eta-\bar{\eta}|^2$ in the vector-valued case. 
In many cases $\bar{\mathcal{E}}$ is assumed to be zero, but we keep it here to be most general.
Here, we note that in contrast to the monoatomic case \eqref{moments}, we have an additional temperature related to the degrees of freedom in internal energy $T_{int}$. From the physical principle equipartition of internal energy, one expects that in equilibrium these two temperatures are the same. To achieve this, there are different strategies in modelling
\begin{itemize}
\item Relaxation to an equilibrium distribution with equal temperature
\end{itemize}
We define the equilibrium temperature as
\begin{align*}
T_{equ}= \frac{d}{d+l} T_{tr} + \frac{l}{d+l} T_{int}
\end{align*}
\note[Review2]{I changed 3 into d in the equilibrium temperature}
Then, we consider the equilibrium distribution
\begin{align}
M_{equ} = \frac{n~ \Lambda_{l}}{\sqrt{2\pi \frac{T_{equ}}{m}}^{d} (T_{equ})^{l/2}} e^{- \frac{m |v-u|^2}{2 T_{equ}} - \frac{e(\mathcal{E})}{T_{equ}}}.
\label{equ_distr}
\end{align}
\note[Review2]{I exchanged $\rho$ into $n$ in (25) and.}
\add[Review2]{with $\Lambda_{l}$ being a constant ensuring that the integral of $M_{equ}$ with respect to v and $\mathcal{E}$ is equal to $n$}.
and the BGK model
\begin{align*}
\partial_t f + v \cdot \nabla_x f = \nu n (M_{equ} -f)
\end{align*}
If we multiply this equation with respect to $e(\mathcal{E})$, we obtain the macroscopic equation in the space-homogeneous case. 
\begin{align*}
\partial_t T_{int} = \frac{\nu m}{l} (T_{equ}- T_{int})= \frac{\nu m}{l} \frac{d}{d+l} (T_{tr}-T_{int})
\end{align*}
so $T_{int}$ relaxes towards $T_{tr}$ with a rate dependent on the collision frequency $\nu$.
\begin{itemize}
\item Relaxation of the temperature due to a convex combination of temperatures in the Maxwell distribution
\end{itemize}
We define
 $T_{rel} = \theta T_{equ} + (1- \theta) T_{int},$ with $0< \theta \leq 1$. Then, we consider the following Maxwell distribution  
$$ \widetilde{G}[f]= \frac{n ~ \Lambda_l}{\sqrt{2 \pi \frac{T}{m}}^d} \frac{1}{\sqrt{ T_{rel}}^l} \exp \left( - \frac{1}{2} \frac{m |v-u|^2}{ T} - \frac{I^{\frac{l}{2}}}{ T_{rel}} \right),$$ 
with the temperature $T = (1- \theta)  T_{tr}  + \theta  T_{equ} $.  $\Lambda_{l}$ is a constant ensuring that the integral of $\widetilde{G}[f]$ with respect to $v$ and $I$ is equal to the density $n$. Then the model is given by
$$ \partial_t f + v \cdot \nabla_x f = A_{\nu} (\widetilde{G}[f] - f)$$
with the collision frequency $A_{\nu}$.
If we choose $f^{(2)}$ as distribution function, this corresponds to the model in \cite{AndriesPerthame2001}.
\change[rephrasing]{For this model one can prove conservation of the number of particles, momentum and total energy, and an entropy inequality such that the equilibrium is characterized by a Maxwell distribution with equal temperatures $T_{equ}= T_{tr} = T_{int}$, for details see section 3 in}{For this model one can show conservation of the number of particles, momentum and total energy. Moreover, one can prove an entropy inequality. Here,  the equilibrium is characterized by a Maxwell distribution with equal temperatures $T_{equ}= T_{tr} = T_{int}$, for details see section 3 in} \cite{AndriesPerthame2001}.\change[rephrasing]{ It is also ensured that there exists a unique mild solution to this model. This is proven in}{ One can also show that there exists a unique mild solution to this model. This is proven in} \cite{Yunexpoly}.

\change[rephrasing]{The convex combination in $\theta$ takes into account that $T_{tr}$ and $T_{int}$ relax towards the common value $T_{equ}$. In the space-homogeneous case we see that we get the following macroscopic equations}{With the convex combination in $\theta$ one takes into account that $T_{tr}$ and $T_{int}$ relax to the common value $T_{equ}$. In the space-homogeneous case one can compute the following macroscopic equations}
\begin{align}
\begin{split}
\label{asd}
\partial_t T_{tr} &= A_{\nu}  ( T_{tr} (1- \theta) + \theta T_{equ} - T_{tr})= A_{\nu} \theta (T_{equ} - T_{tr}),
\\
\partial_t T_{int} &= A_{\nu} \theta (T_{equ} - T_{int}) .
\end{split}
\end{align}
 These macroscopic equations describe a relaxation of $T_{tr}$ and $T_{int}$ towards $T_{equ}$ with a speed depending on the additional parameter $\theta$. 
\add[Review2]{  We see that the model captures the regime where this relaxation of the temperatures is slower than the relaxation of the distribution function to a Maxwell distribution since $\theta$ satisfies $\theta \leq 1$, so it reduces the speed of relaxation from $A_{\nu}$ to $A_{\nu}  \theta$. 
}
This model satisfies the following assymptotic behaviour proven in \cite{Yunpoly} in the space-homogeneous case.

\begin{thm}
Let $0<\theta \leq 1.$ The distribution function for the spatially homogeneous case converges to equilibrium with the following rate:
$$||f(t) - M_{equ}||_{L^1(dv dI)} \leq e^{- \frac{\theta}{2} A_{\nu} t }\sqrt{2 H(f_0 |M_{equ})},$$
with the relative entropy $H(f|g)= \int \int f \ln \frac{f}{g} dv dI$ for two functions $f$ and $g$, and the Maxwell distribution $M_{equ}$ given by \eqref{equ_distr}.
\end{thm}
\begin{itemize}
\item Relaxation of the temperatures with an additional kinetic equation
\end{itemize}
This concept was introduced in \cite{Bernard} for the distribution function $f^{(4)}$. Here, we discribe the time evolution in the following way
\begin{align} \begin{split} \label{BGKone}
\partial_t f + v\cdot\nabla_x   f   = \nu n (M[f] - f) 
\end{split}
\end{align}
with the Maxwell distribution
\begin{align} 
\begin{split}
M(x,v,\mathcal{E},t) = \frac{n}{\sqrt{2 \pi \frac{\Lambda}{m}}^d } \frac{1}{\sqrt{2 \pi \frac{\Theta}{m}}^{l}} \exp({- \frac{|v-u|^2}{2 \frac{\Lambda}{m}}}- \frac{e(\mathcal{E})}{2 \frac{\Theta}{m}}), 
\end{split}
\label{BGKmixone}
\end{align}
where $\nu n$ is the collision frequency. Here, there appear to additional artificial temperatures $\Lambda$ and $\Theta$. In order to describe the time evolution of these two temperatures, 
we couple this kinetic equation with an algebraic equation for conservation of internal energy 
\begin{align}
\frac{d}{2} n \Lambda = \frac{d}{2} n T_{tr} +\frac{l}{2} n T_{int}- \frac{l}{2} n\Theta,  \label{internalone}
\end{align}  
and a relaxation equation ensuring that the two temperatures $\Lambda$ and $\Theta$ relax to the same value in equilibrium
 \begin{align}
 \begin{split}
 \partial_t M + v \cdot \nabla_x M &= \frac{\nu n}{Z_r} \frac{d+l}{d} (M_{equ} - M) \\ 
 \Theta(0)&= \Theta^0
 \end{split} 
 \label{kin_Tempone}
 \end{align}
 where $Z_r$ is a given parameter corresponding to the different rates of decay of translational and rotational/vibrational degrees of freedom. 
 Here $M$ is given by
\begin{align} 
M(x,v,\mathcal{E},t) = \frac{n}{\sqrt{2 \pi \frac{\Lambda}{m}}^d} \frac{1}{\sqrt{2 \pi \frac{\Theta}{m}}^{l}} \exp \left({- \frac{|v-u|^2}{2 \frac{\Lambda}{m}}}- \frac{e(\mathcal{E})}{2 \frac{\Theta}{m}}\right), 
\label{Max_equone}
\end{align}
Note that we have
\begin{align}
T_{equ}= \frac{d \Lambda+ l \Theta}{d+l}= \frac{d T_{tr} + l T_{int}}{d+l}.
\label{equ_tempone}
\end{align}
The second equality follows from \eqref{internalone}. The equation \eqref{kin_Tempone} is used to involve the temperature $\Theta$. If we multiply \eqref{kin_Tempone} by $e(\mathcal{E})$, integrate with respect to $v$ and $\mathcal{E}$ and use \eqref{equ_tempone}, we obtain  
\begin{align}
\begin{split}
\partial_t(n \Theta) +   \nabla_x\cdot (n \Theta u) = \frac{\nu n}{Z_r} n (\Lambda - \Theta) 
\end{split}
\label{relaxone}
\end{align} 
a relaxation to a common value with a speed $Z_r$ not restricted to a slower speed as in \eqref{asd}. \add[Review2]{ Therefore, the relaxation of the two temperatures to a common value can be slower or faster than the relaxation of the distribution function to the Maxwell distribution. This depends on the choice of $Z_r$.}
The initial data of $\Lambda$ is determined using \eqref{internalone}.\change[rephrasing]{ We see that in this model the term $\frac{\nu n}{Z_r} \frac{d+l}{d} (M_{equ} - M)$ describes the relaxation of the two temperatures $\Lambda$ and $\Theta$ to a common value. So the effect of the relaxation to common temperatures here is done by coupling the BGK equation with an additional kinetic equation. 
If we choose $f^{(4)}$ as distribution function, this corresponds to the model in}{We see that in this model the term $\frac{\nu n}{Z_r} \frac{d+l}{d} (M_{equ} - M)$ plays the role to produce the relaxation of the two temperatures $\Lambda$ and $\Theta$ to the same value. So in this model the effect of the relaxation to equal temperatures  is done by coupling the BGK equation with an additional kinetic equation. 
If we choose $f^{(4)}$ as distribution function, this corresponds to the model in} \cite{Bernard}.
\change[rephrasing]{For this model one can prove conservation of the number of particles, momentum and total energy, and an entropy inequality such that the equilibrium is characterized by a Maxwell distribution with equal temperatures $T_{equ}= T_{tr} = T_{int}$, for details see }{For this model one can prove conservation of the number of particles, momentum and total energy, and also an entropy inequality. Additionally, the equilibrium can be characterized by a Maxwell distribution with equal temperatures $T_{equ}= T_{tr} = T_{int}$, for details see } \cite{Bernard, Pirner_poly_con}. \change[rephrasing]{The existence of mild solution can be proven similar to the existence in the momatomic case}{The existence of a unique mild solution can be proven similar to the existence in the momatomic case} \cite{Perthame}.

This model satisfies the following asymptotic behaviour in the space-homogeneous case proven in \cite{Pirner_poly_con} for $f=f^{(4)}.$
\begin{thm}
\label{conv3}
 Assume that $(f^{(4)}, M)$ is a solution of \eqref{BGKone} coupled with \eqref{kin_Tempone} and \eqref{internalone}. Then, in the space homogeneous case, we have the following convergence rate of the distribution functions $f$:
 {\small
\begin{align*}
||f^{(4)} - M_{equ}||_{L^1(dv d\eta)}  \leq 4 e^{-\frac{1}{4} \widetilde{C} t} \left( H(f^{(4)}_0|M^0_{equ})+ 2 \max \lbrace 1, z\rbrace H(M^0|M^0_{equ}) \right)^{\frac{1}{2}}.
\end{align*}}
where $\widetilde{C}$ is given by $$\widetilde{C}= \min \left\lbrace  \nu n^{(4)}  ,  \frac{\nu n^{(4)}}{Z_r} \frac{d+l}{d}  \right\rbrace,$$
and the index $0$ denotes the value at $t=0$.
\end{thm}
\begin{itemize}
\item Relaxation of the temperatures with an additional relaxation term
\end{itemize}
This concept was introduced in \cite{Blaga}. Here, we add an additional relaxation term into the right-hand side
\begin{align*}
\partial_t f + v \cdot \nabla_x f = \frac{1}{\tau} (m_2 -f) + \frac{1}{Z_{\nu} \tau} (M_{equ} - m_2)
\end{align*}
where $\tau$ is the relaxation time of $f$ towards a Maxwell distribution with the temperatures $T_{tr}$ and $T_{int}$ given by
\begin{align*}
m_2 = \frac{n}{\sqrt{2 \pi \frac{T_{tr}}{m}}^d} \frac{1}{\sqrt{2 \pi \frac{T_{int}}{m}}^{l}} \exp({- \frac{|v-u|^2}{2 \frac{T_{tr}}{m}}}- \frac{e(\mathcal{E})}{2 \frac{T_{int}}{m}}), 
\end{align*}
and $\tau Z_{\nu}$ with $Z_{\nu}>1$ the relaxation time of this Maxwell distribution to the equilibrium distribution $M_{equ}$  with equal temperatures given by \eqref{Max_equone}.  So the relaxation to equilibrium is divided into two parts, first a relaxation towards an intermediate equilibrium distribution where the temperatures $T_{tr}$ and $T_{int}$ are separate, then to the final equilibrium with common temperatures.
\subsection{Summary of existing BGK models for gas mixtures of polyatomic molecules in the literature}
Now, for the gas mixture case, we will present different models \cite{Bisi_poly_mix,Pirner_poly_con,Blaga} combining different ansatzes from the one species polyatomic  case and the mixture modelling. 
\subsubsection{A BGK model for mixtures of polyatomic gases with one relaxation term}
The BGK model we consider in this section was introduced by Bisi, Monaco and Soares in \cite{Bisi_poly_mix}. We introduce here two distribution functions with scalar continuous dependency on the degrees of freedom of internal energy $f_1(x,v,t,I)$ and $f_2(x,v,t,I)$. Then the time evolution of these distribution functions is described by two kinetic equations with one relaxation term on the right-hand side to the equilibrium distribution with common temperatures
\begin{align*}
\partial_t f_k + v \cdot \nabla_x f_k = \nu_k (M^k-f_k), \quad k=1,2
\end{align*}
with the Maxwell distributions 
\begin{align*}
M^k(v,I)= \frac{\tilde{n}_k}{q^k(\tilde{T})} \left( \frac{m_k}{2 \pi \tilde{T}} \right)^{\frac{d}{2}} \exp \left(-\frac{1}{\tilde{T}} ( \frac{m_k}{2} |v- \tilde{u}|^2 + I)\right), \quad k=1,2
\end{align*}
with the partition function $q^k(\tilde{T}) = \int_0^{\infty} \phi^k(I) \exp(-\frac{I}{T}) dI$. Then, the parameters $\tilde{n}_k, \tilde{u}$ and $\tilde{T}$ will be determined to have conservation of mass, total momentum and total energy. For the computation and the detailed expression see \cite{Bisi_poly_mix}. For this model also an entropy inequality can be proven in the space-homogeneous case, see \cite{Bisi_poly_mix}. \add[Review2]{Transport coefficients in the hydrodynamical limit of this model can be find in section 5 of} \cite{Bisi_poly_mix}.
\subsubsection{A BGK model for mixtures of polyatomic gases  with two relaxation terms}
In this section, we present the model developed in \cite{Pirner_poly_con}. \change[rephrasing]{This model has a vector-valued dependency on the internal energy. For this we introduce the total number of different rotational and vibrational degrees of freedom $M$; and $l_k$ the number of internal degrees of freedom of species $k,~k=1,2$. Further $\eta \in \mathbb{R}^M$ is the variable for the internal energy degrees of freedom, whereas $\eta_{l_k} \in \mathbb{R}^M$ coincides with $\eta$ in the components corresponding to the internal degrees of freedom of species $k$ and is zero in the other components.}{ This model has a vector-valued dependency on the internal energy. For this we introduce two numbers related to the degrees of freedom in internal energy. One is the total number of different rotational and vibrational degrees of freedom $M$; and the other is $l_k$, the number of internal degrees of freedom of species $k,~k=1,2$. Moreover, $\eta \in \mathbb{R}^M$ is the variable for the internal energy degrees of freedom, whereas $\eta_{l_k} \in \mathbb{R}^M$ coincides with $\eta$ in the components corresponding to the internal degrees of freedom of species $k$ and is zero in all the other components.} In this way, it is possible that the two species can have a different number of degrees of freedom in internal energy. Then, we have distribution functions $f_1(x,v,t,\eta_{l_1})$ and $f_2(x,v,t, \eta_{l_2})$. Their time evolution is described by
\begin{align} \begin{split} \label{BGK3}
\partial_t f_1 + v\cdot\nabla_x   f_1   &= \nu_{11} n_1 (M_1 - f_1) + \nu_{12} n_2 (M_{12}- f_1),
\\ 
\partial_t f_2 + v\cdot\nabla_x   f_2 &=\nu_{22} n_2 (M_2 - f_2) + \nu_{21} n_1 (M_{21}- f_2), \\
\end{split}
\end{align}
with the Maxwell distributions
\begin{align} 
\begin{split}
M_k(x,v,\eta_{l_k},t) &= \frac{n_k}{\sqrt{2 \pi \frac{\Lambda_k}{m_k}}^d } \frac{1}{\sqrt{2 \pi \frac{\Theta_k}{m_k}}^{l_k}} \exp({- \frac{|v-u_k|^2}{2 \frac{\Lambda_k}{m_k}}}- \frac{|\eta_{l_k}- \bar{\eta}_{l_k}|^2}{2 \frac{\Theta_k}{m_k}}), 
\\
M_{kj}(x,v,\eta_{l_k},t) &= \frac{n_{kj}}{\sqrt{2 \pi \frac{\Lambda_{kj}}{m_k}}^d } \frac{1}{\sqrt{2 \pi \frac{\Theta_{kj}}{m_k}}^{l_k}} \exp({- \frac{|v-u_{kj}|^2}{2 \frac{\Lambda_{kj}}{m_k}}}- \frac{|\eta_{l_k}- \bar{\eta}_{l_k,kj}|^2}{2 \frac{\Theta_{kj}}{m_k}}), 
\end{split}
\label{BGKmix3}
\end{align}
for $ j,k =1,2, j \neq k$ with the conditions 
\begin{equation} 
\nu_{12}=\varepsilon \nu_{21}, \quad 0 < \frac{l_1}{l_1+l_2}\varepsilon \leq 1.
\label{coll3}
\end{equation}
The equation is coupled with conservation of internal energy 
\eqref{internalone} for each species, and an additional relaxation equation 
 \begin{align}
 \begin{split}
 \partial_t M_k + v \cdot \nabla_x M_k = \frac{\nu_{kk} n_k}{Z_r^k} \frac{d+l_k}{d} (M_{equ,k} - M_k)+ \nu_{kj} n_j  (\widetilde{M}_{kj} - M_k) , \\
 \Theta_k(0)= \Theta_k^0
 \end{split} 
 \label{kin_Temp3}
 \end{align}
for $j,k=1,2, j \neq k$.
 $M_{equ,k}$ is given by \eqref{Max_equone} for each species.  The additional $\widetilde{M}_{kj}$ is defined by
\begin{align}
\widetilde{M}_{kj}= \frac{n_k}{\sqrt{2 \pi \frac{T_{kj}}{m_k}}^{d+l_k}} \exp \left(- \frac{m_k |v-u_{kj}|^2}{2 T_{kj}}- \frac{m_k|\eta_{l_k}- \bar{\eta}_{l_k,kj}|^2}{2 T_{kj}} \right), \quad k=1,2.
\label{Max_equ3}
\end{align}
where $T_{kj}$ is  given by 
\begin{align}
T_{kj}:= \frac{d \Lambda_{kj} + l_k \Theta_{kj}}{d+l_k}.
\label{equ_temp3}
\end{align}

For a certain choice of $\Lambda_{kj}$ and $\Theta_{kj}$ one can prove conservation of mass, total momentum and total  energy. For details see \cite{Pirner_poly_con}. The existence of solutions for this model can be proven in the same way as it is done in \cite{Pirner} for the monoatomic case. In \cite{Pirner_poly_con} they also prove an entropy inequality and the following decay to equilibrium

\begin{thm}
\label{conv3}
 Assume that $(f_1, f_2, M_1, M_2)$ is a solution of \eqref{BGK3} coupled with \eqref{kin_Temp3} and \eqref{internalone}. Then, in the space homogeneous case, we have the following convergence rate of the distribution functions $f_1$ and $f_2$:
 {\small
\begin{align*}
||f_k - \widetilde{M}_k||_{L^1(dv d\eta_{l_k})}  \leq 4 e^{-\frac{1}{4} \widetilde{C} t} \left(\sum_{k=1}^2 \left( H_k(f_k^0|\widetilde{M}_k^0)+ 2 \max \lbrace 1, z_1, z_2 \rbrace H_k(M_k^0|\widetilde{M}_k^0)\right) \right)^{\frac{1}{2}}.
\end{align*}}
where $\widetilde{C}$ is given by $$\widetilde{C}= \min \left\lbrace  \nu_{11} n_1 +  \nu_{12} n_2 , \nu_{22} n_2 + \nu_{21} n_1, \frac{\nu_{11} n_1}{z_1} + \nu_{12} n_2, \frac{\nu_{22} n_2}{z_2} + \nu_{21} n_1 \right\rbrace,$$
and the index $0$ denotes the value at time $t=0$.
\end{thm}
There are also numerical results for this model in \cite{Karlsruhe}.

\subsection{BGK model for mixtures of polyatomic gases with intermediate relaxation terms}
The model in \cite{Blaga} extends the idea of additional relaxation terms with intermediate equilibrium distributions from the one species case to gas mixtures. The model is of the form
\begin{align*}
\partial_t f_k + v \cdot \partial_x f_k = \frac{1}{\tau} (m_{s_1} - f_s) + \frac{1}{Z_r \tau} (m_{s_2} -m_{s_1}) + \frac{1}{Z_{\nu} \tau} (\tilde{M}_k -m_{s_2}), \quad k=1,2
\end{align*}
with $Z_r,Z_{\nu} >1$ and intermediate equilibrium distributions $m_{s_1}$ and $m_{s_2}$.
The detailed expressions of the intermediate equilibrium distributions can be found in \cite{Blaga} with a proof of the conservation properties. With standard methods one can also prove an entropy inequality. \add[Review2]{Transport coefficients of the hydrodynamic regime of this model can be found in section 4 of} \cite{Blaga}.

\section{Conclusions}

\add[Review1]{This paper reviews various existing BGK models for gas mixtures of mono and polyatomic molecules from the literature. In the case of monoatomic particles   two types of models are presented. 

One contains only one relaxation term on the right-hand side taking into account all types of interactions in one relaxation term. The other type of model separates the inter- and intra-species interactions by writing a sum of relaxation terms on the right-hand side.  For both types of models a review on theoretical results concerning existence of solutions and convergence to equilibrium are given. The results on convergence to equilibrium consider both the space-homogeneous case for the full non-linear model and the space-inhomogeneous case for a linearized model.

 In the polyatomic case first different ansatzes of modelling the degrees of freedom in internal energy for one species are considered: discrete or continuous; and scalar- or vector-valued.

  Next, different ansatzes for modelling the relaxation of the temperature related to translational degrees of freedom and the temperature related to internal degrees of freedom to the same value for one species are presented. 
Here, also theoretical results on the convergence to equilibrium are presented. 

Finally, three existing BGK  models in the literature concerning gas mixtures of polyatomic molecules using the different ansatzes of polyatomic and gas mixture modelling are presented with existing theoretical results on the convergence to equilibrium.

 However, BGK- type models often lack on correct parameters in the continuum limit like the Prandtl number. Therefore these models can be used as a basis for more extended models like ES-BGK models or velocity dependent collision frequency. As a future work the Chapman- Enskog expansion of the missing models can be computed and then the transport coefficients of all these models can be compared and eventually extended to match all parameters in the macroscopic equations. Here, the free parameters in the BGK model for monoatomic molecules with a sum of interaction terms might be useful.  }
\vspace{6pt} 

%
%

{\bf funding:} This research was funded by the Alexander von Humboldt foundation.

\end{document}